\documentclass[12pt]{article}

\usepackage{amssymb,amsfonts,amsmath}

\newtheorem{theorem}{Theorem}

\newtheorem{remark}{Remark}

\newcommand{\rom}[1]{{\rm #1}}

\newcommand{\Ffin}{\mathcal F_{\mathrm {fin}}}

\newcommand{\D}{{\mathcal D}}
\newcommand{\N}{{\mathbb N}}

\newcommand{\Z}{{\mathbb Z}}
\newcommand{\R}{{\mathbb R}}

\newcommand{\la}{\langle}
\newcommand{\ra}{\rangle}

\newcommand{\hotimes}{\hat\otimes}

\newcommand{\ZZ}{\Z_{+,\,0}^\infty}

\begin{document}

\begin{center}

{\Large \bf An equivalent representation \\ of the Jacobi field of 
a L\'evy process }\\[10mm]
{\large\bf E. Lytvynov}\\[2mm]
 Department of Mathematics\\
University of Wales Swansea \\
Singleton Park\\
Swansea SA2 8PP\\
U.K.\\
E-mail: e.lytvynov@swansea.ac.uk

\end{center}

$\text{}$\hfill {\it Dedicated to  Professor Yuri Berezansky}\\
$\text{}$\hfill {\it on his 80th birhday}

\begin{abstract}

In \cite{BLM}, the Jacobi field of a L\'evy process was derived. 
This  field consists of commuting self-adjoint operators  acting in an
extended (interacting) Fock space. However, these operators have a quite complicated structure. In this note, using  ideas from 
\cite{AFS, LSWN}, we obtain a unitary equivalent representation of the Jacobi field of a L\'evy process. In this representation, the operators
act in a usual symmetric Fock space and have a much simpler structure.

\end{abstract}

\noindent {\it AMS Mathematics Subject Classification}:
 60G20, 60G51,   60H40, 47B36


\section{L\'evy process and its Jacobi field}

The notion of a Jacobi field in the Fock space first appeared in the works by
Berezansky and Koshmanenko \cite{BerKosh1,BerKosh2},
devoted to the axiomatic quantum field theory, and then was further developed by
Br\"uning  (see e.g.\ \cite{Bru}).
These works, however, did not contain any relations with probability measures.  A detailed study of a general commutative Jacobi  field in the Fock space and a corresponding spectral  measure was carried out in a serious of works by Berezansky, see e.g.\ \cite{bere, berre} and the references therein.

In \cite{BLM} (see also \cite{L, Lproc}), the Jacobi field of a L\'evy process on a general manifold $X$ was studied. Let us shortly recall
these results. 

Let $X$ be a complete, connected, oriented $C^\infty$
(non-compact) Riemannian manifold and let ${\cal B}(X)$ be the
Borel $\sigma$-algebra on $X$. Let $\sigma$ be a Radon measure on
$(X,{\cal B}(X))$ that is non-atomic and non-degenerate (i.e.,
$\sigma(O)>0$ for any open set $O\subset X$). As a typical example
of  measure $\sigma$, one can take the volume measure on $X$.

We denote by $\D$
the space $C_0^\infty(X)$ of all infinitely differentiable,
real-valued functions on $X$ with compact support. It is known
that $\D$ can be endowed with a topology of a nuclear space. Thus, we can consider the standard
nuclear triple $$ \D\subset L^2(X,\sigma)\subset\D',$$  where $\D'$
is the dual space of $\D$ with respect to the zero space
$L^2(X,\sigma)$. (Here and below, all the linear spaces we deal with are real.) The dual pairing between $\omega\in\D'$ and
$\varphi\in\D$  will be denoted by $\la\omega,\varphi\ra$.  
We  denote the cylinder $\sigma$-algebra on $\D'$ by ${\cal C}(\D')$.

Let $\nu$ be a measure on $(\R,{\cal B}(\R))$ whose support contains an infinite number of points and assume $\nu(\{0\})=0$. 
Let $$\tilde\nu(ds):=s^2\,\nu(ds).$$ We further assume  that
$\tilde \nu$ is a finite measure on $(\R,{\cal B}(\R))$, and
moreover,   there exists an $\varepsilon>0$
 such that
\begin{equation}\label{4335ew4} \int_{\R} \exp\big(\varepsilon
|s|\big)\,\tilde\nu(ds)<\infty.\end{equation}   

We now
define a centered L\'evy process on $X$  (without Gaussian part)
as a generalized process
on ${\cal D}'$ whose law is the probability measure
$\mu$ on $({\cal D}',{\cal C}({\cal D}'))$ given
by its Fourier transform
 \begin{equation}\label{rew4w}
\int_{{\cal D}'} e^{i\la \omega,\varphi\ra}\,
\mu(d\omega)=\exp\bigg[\int_{\R\times X}(e^{is\varphi(x)}-1-is\varphi(x))\,\nu(ds)\,\sigma(dx)\bigg],\qquad
\varphi\in {\cal D}.\end{equation} 

Thus, $\nu$ is the L\'evy measure of the L\'evy process $\mu$.
Without loss of generality, we can suppose that $\tilde \nu$ is a probability measure on $\R$. (Indeed, if this is not the case, define
$\nu':=c^{-1}\nu$ and $\sigma':=c\sigma$, where $c:=\tilde\nu(\R)$.)

It follows from \eqref{4335ew4} that the measure $\tilde\nu$
has all moments finite, and furthermore, the set of all polynomials is dense in 
$L^2(\R,\tilde\nu)$. Therefore, by virtue of \cite{b}, there exists 
a unique (infinite) Jacobi matrix $$ J=\left(\begin{matrix} a_0& b_1& 0& 0&0&\dots\\
b_1&a_1& b_2&0&0&\dots\\ 0&b_2&a_2&b_3&0&\dots\\
0&0& b_3&a_3&b_4&\dots\\ \vdots&\vdots& \vdots&\vdots&\vdots&\ddots
\end{matrix}\right),\qquad a_n\in\R,\ b_n>0,$$ whose spectral measure is $\tilde\nu$.

Next, we denote by ${\cal P}({\cal D}')$ the set of continuous
polynomials on ${\cal D}'$, i.e., functions on ${\cal D}'$ of the
form $F(\omega)=\sum_{i=0}^n\la\omega^{\otimes i},f_{i}\ra$, $
\omega^{\otimes 0}{:=}1,\ f_{i}\in{\cal D}^{\hotimes i}$,
$i=0,\dots,n$, $n\in\Z_+$.  Here, $\hat\otimes$ stands for
symmetric tensor product.  The greatest number $i$ for which
$f^{(i)}\ne0$ is called the power of a polynomial. We denote by
${\cal P}_n({\cal D}')$ the set of continuous polynomials of power
$\le n$.

By \eqref{4335ew4}, \eqref{rew4w}, and \cite[Sect.~11]{Sko},
${\cal P}({\cal D}')$ is a dense subset of $ L^2({\cal
D}',\mu)$. Let ${\cal P}^\sim _n({\cal D}')$ denote the closure of ${\cal
P}_n({\cal D}')$ in $L^2({\cal D}',\mu)$, let
${\bf P}_n({\cal D}')$, $n\in\N$, denote the orthogonal difference
${\cal P }^\sim_n({\cal D}')\ominus{\cal P}^\sim_{n-1}({\cal
D}')$, and let ${\bf P}_0({\cal D}'){:=}{\cal P }^\sim_0({\cal
D}')$. Then, we evidently have: \begin{equation}\label{zgzguzug}
L^2({\cal D}',\mu)=\bigoplus_{n=0}^\infty{\bf
P}_n({\cal D}'). \end{equation}

The set of all projections ${:}\la\cdot^{\otimes n},f_n\ra{:}$ of
continuous monomials $\la \cdot^{\otimes n},f_n\ra$,
$f_n\in\D^{\hat\otimes n}$, onto ${\bf P}_n({\cal D }')$ is dense
in ${\bf P}_n({\cal D}')$. For each $n\in\N$, we define a Hilbert
space ${\frak F}_n$
 as the closure of the set ${\cal D}^{\hotimes n}$ in the
 norm generated by the scalar product \begin{equation}\label{gffdz}
 (
 f_n,g_n)_{{\frak F}_n}{:=}\frac1{n!}\, \int_{{\cal D}'}{:}\la \omega^{\otimes n},f_n\ra{:}
 \, {:}\la\omega^{\otimes n},g_n\ra{:}\,\mu(d\omega),\qquad f_n,g_n\in{\cal D}^{\hotimes n}.
 \end{equation} Denote  \begin{equation}\label{qztztzt}{\frak
F}{:=}\bigoplus_{n=0}^\infty{\frak F}_n\,n!,\end{equation} where
${\frak F}_0{:=}\R$. By \eqref{zgzguzug}--\eqref{qztztzt}, we get
 the unitary operator $${\cal U}:{\frak F}\to
L^2(\D',\mu)$$  that is defined through $ {\cal
U}f_n{:=}{:}\la\cdot^{\otimes n},f_n\ra{:}$, $f_n\in{\cal
D}^{\hotimes n}$, $n\in\Z_+$, and then extended by linearity and
continuity to the whole space ${\frak F}$\rom.

An explicit formula for the scalar product $(\cdot,\cdot)_{{\frak F}_n}$ 
looks as follows. We denote by $\ZZ $ the set of all sequences $\alpha$ of the form
$$\alpha=(\alpha_1,\alpha_2,\dots,\alpha_n,0,0,\dots),\qquad
\alpha_i\in\Z_+,\ n\in\N.$$ Let $|\alpha|{:=}\sum_{i=1}^\infty
\alpha_i$. For each $\alpha\in\ZZ$, $1\alpha_1+2\alpha_2+\dots=n$,
$n\in\N$, and for any function $f_n:X^n\to\R$ we define a function
$D_\alpha f_n:X^{|\alpha|}\to\R$ by setting \begin{align*}(D_\alpha
f_n)(x_1,\dots,x_{|\alpha|}){:=}& f(x_1,\dots,x_{\alpha_1},
\underbrace{x_{\alpha_1+1},x_{\alpha_1+1}}_{\text{2 times }},
\underbrace{x_{\alpha_1+2},x_{\alpha_1+2}}_{\text{2 times
}},\dots,
\underbrace{x_{\alpha_1+\alpha_2},x_{\alpha_1+\alpha_2}}_{\text{2
times }},\notag\\ &\quad
\underbrace{x_{\alpha_1+\alpha_2+1},x_{\alpha_1+\alpha_2+1},x_{\alpha_1+\alpha_2+1}}_{\text{3
times }},\dots).\end{align*}

We have (cf.\ \cite{L}):

\begin{theorem}\label{oifvoh}
 For any
$f^{(n)},g^{(n)}\in{\cal D}^{\hat\otimes n}$\rom, we have\rom:
\begin{align*} (f^{(n)},g^{(n)})_{{\frak F}_n}& =\sum_{\alpha\in\ZZ:\,
1\alpha_1+2\alpha_2+\dots=n}K_{\alpha}
\int_{X^{|\alpha|}}(D_\alpha f_n)(x_1,\dots,x_{|\alpha|})\notag
\\ & \quad \times (D_\alpha g_n)(x_1,\dots,x_{|\alpha|}) \,\sigma^{\otimes
|\alpha|}(dx_1,\dots,dx_{|\alpha|}),\end{align*} where
\begin{equation}\label{iugu} K_{\alpha}= \frac{(1\alpha_1+2\alpha_2+\dotsm)!}
{\alpha_1!\,\alpha_2!\dotsm}\,\prod_{k\ge2}\bigg(\frac{\prod_{i=1}^{k-1}
b_i}{k!}\bigg)^{2\alpha_k}. \end{equation}
\end{theorem}

Next, we find the elements which belong to the space ${\frak F}_n$
after the completion of ${\cal D}^{\hat\otimes n}$. To this end, 
we define, for each $\alpha\in\ZZ$\rom, the Hilbert
space $$ L^2_\alpha(X^{|\alpha|},\sigma^{\otimes
|\alpha|}){:=}L^2(X,\sigma)^{\hotimes \alpha_1}\otimes
L^2(X,\sigma)^{\hotimes\alpha_2}\otimes\dotsm\,.$$ Define a
mapping $$ U^{(n)}: {\cal D}^{\hat\otimes n}\to
\bigoplus_{\alpha\in\ZZ:\,
1\alpha_1+2\alpha_2+\dots=n}L^2_\alpha(X^{|\alpha|},\sigma^{\otimes|\alpha|})K_{\alpha}$$
by setting\rom, for each $f^{(n)}\in{\cal D}^{\hat\otimes n}$\rom, the
$L^2_\alpha(X^{|\alpha|},\sigma^{\otimes|\alpha|})K_\alpha$-coordinate of
$U^{(n)}f^{(n)}$ to be $D_\alpha f^{(n)}$.
By virtue of  Theorem~\ref{oifvoh}, $U^{(n)}$ may be extended
by continuity to an isometric mapping of ${\frak F}_n$ into $$\bigoplus_{\alpha\in\ZZ:\,
1\alpha_1+2\alpha_2+\dots=n}L^2_\alpha(X^{|\alpha|},\sigma^{\otimes|\alpha|})K_{\alpha}.$$
Furthermore, we have (cf.\ \cite{beme1,L}):

\begin{theorem}\label{ojivfed}
The mapping $$U^{(n)}:{\frak F}_n\to\bigoplus_{\alpha\in\ZZ:\,
1\alpha_1+2\alpha_2+\dots=n}L^2_\alpha(X^{|\alpha|},\sigma^{\otimes|\alpha|})K_{\alpha}$$
 is  a unitary opertator\rom.
\end{theorem}

By virtue of Theorem~\ref{ojivfed} and \eqref{qztztzt}, we can identify ${\frak F}_n$ with the space $$\bigoplus_{\alpha\in\ZZ:\,
1\alpha_1+2\alpha_2+\dots=n}L^2_\alpha(X^{|\alpha|},\sigma^{\otimes|\alpha|})K_{\alpha}$$ and the space $\frak F$ with 
$$ \bigoplus_{\alpha\in\ZZ}L^2_\alpha(X^{|\alpha|},\sigma^{\otimes|\alpha|})K_{\alpha}
(1\alpha_1+2\alpha_2+\cdots)!\,. $$ For a vector $f\in{\frak F}$, we will denote its  $\alpha$-coordinate   by $f_\alpha$ . 

Note that, for for $\alpha=(n,0,0,\dots)$, we have
$$L^2_\alpha(X^{|\alpha|},\sigma^{\otimes|\alpha|})=L^2(X,\sigma)^{\hat\otimes n},\quad K_\alpha=1,\quad (1\alpha_1+2\alpha_2+\cdots)!=n!.$$
Hence, the space $\frak F$ contains the symmetric Fock space
$${\cal F}(L^2(X,\sigma))=\bigoplus_{n=0}^\infty L^2(X,\sigma)^{\hat\otimes n} n!$$ as a proper subspace.
Therefore, we call $\frak F$ an extended
Fock space. We also note that the space $\frak F$ satisfies the axioms of an interacting Fock space, see \cite{HW}. 

In the space $L^2(\D',\mu)$, we consider, 
for each  $\varphi\in\D$,  the operator $M(\varphi)$ of
multiplication by the function $\la\cdot,\varphi\ra$. Let
$J(\varphi):={\cal U}M(\varphi){\cal U}^{-1}$. Denote by
$\Ffin(\D)$ the set of all vectors of the form
$(f_0,f_1,\dots,f_n,0,0,\dots)$, $f_i\in\D^{\hotimes i}$,
$i=0,\dots,n$, $n\in\Z_+$. Evidently, $\Ffin(\D)$ is a dense
subset of $\frak F$. We have the following theorem, see \cite{BLM}.

\begin{theorem}\label{zgztpopo}
For any $\varphi\in\D$\rom, we have\rom:
\begin{equation}\label{7t667}
\Ffin(\D)\subset\operatorname{Dom}(J(\varphi)),\qquad
J(\varphi)\restriction\Ffin(\D)=J^+(\varphi)+J^0(\varphi)+J^-(\varphi).\end{equation}
Here, $J^+(\varphi)$ is the usual creation operator\rom :
\begin{equation}\label{hbgg} J^+(\varphi) f_n=\varphi\hotimes
f_n,\qquad f_n\in\D^{\hotimes n},\ n\in\Z_+.\end{equation} Next, for each
$f^{(n)}\in\D^{\hat\otimes n}$, $J^0(\varphi)f^{(n)}\in {\frak F}_n$ and
\begin{gather}
(J^0(\xi)f^{(n)})_\alpha(x_1,\dots,x_{|\alpha|})\notag\\=\sum_{k=1}^\infty
\alpha_k a_{k-1}
S_\alpha\big(\xi(x_{\alpha_1+\dots+\alpha_k})(D_\alpha
f^{(n)})(x_1,\dots,x_{|\alpha|})\big)\notag\\
\text{$\sigma^{\otimes|\alpha|}$-a.e.,}\ \alpha\in\ZZ,\
1\alpha_1+2\alpha_2+\dots=n, \label{gfg}\end{gather}
$J^-(\xi)f^{(n)}=0$ if $n=0$\rom, $J^-(\xi)f^{(n)}\in {\frak F}_{n-1}$ if $n\in\N$ and
\begin{gather}
(J^-(\xi)f^{(n)})_\alpha(x_1,\dots,x_{|\alpha|})\notag\\ =n
S_\alpha\bigg(\int_X
\xi(x)(D_{\alpha+1_1}f^{(n)})(x,x_1,\dots,x_{|\alpha|})\,\sigma(dx)\bigg)\notag\\
\text{}+\sum_{k\ge 2 }\frac nk\, \alpha_{k-1}b^2_{k-1} S_\alpha
\big(
\xi(x_{\alpha_1+\dots+\alpha_k})(D_{\alpha-1_{k-1}+1_k}f^{(n)})(x_1,\dots,x_{|\alpha|})\big)\notag\\
\text{$\sigma^{\otimes|\alpha|}$-a.e.},\ \alpha\in\ZZ,\
1\alpha_1+2\alpha_2+\dots=n-1. \label{ggztr}\end{gather} In
formulas \eqref{gfg} and \eqref{ggztr}\rom,  we denoted by
$S_\alpha$  the orthogonal projection of\linebreak
$L^2(X^{|\alpha|},\sigma^{\otimes|\alpha|})$ onto
$L_\alpha^2(X^{|\alpha|},\sigma^{\otimes|\alpha|})$\rom,
$$\alpha\pm
1_n{:=}(\alpha_1,\dots,\alpha_{n-1},\alpha_n\pm1,\alpha_{n+1},\dots),\qquad
\alpha\in\ZZ,\  n\in\N.$$ 

 Finally\rom, each operator $J(\varphi)$\rom, $\varphi\in\D$\rom,
is essentially self-adjoint on $\Ffin(\D)$\rom. 

\end{theorem}

By \eqref{7t667}, the operator $J(\varphi)\restriction\Ffin(\D)$
is a sum of creation, neutral, and  annihilation  operators, and
hence $J(\varphi)\restriction\Ffin(\D)$ has a Jacobi operator's
structure. The family of operators $(J(\varphi))_{\varphi\in \D}$
is called the Jacobi field corresponding to the L\'evy process
$\mu$.

\section{An equivalent representation}

As shown in \cite{berepascal,KL,Ly3,L}, in some cases, the formulas
describing the operators $J^0(\varphi)$ and  $J^-(\varphi)$ can be significantly simplified. However, in the case of a general L\'evy process this is not possible, see \cite{BLM}. We will now present a unitarily equivalent description of the Jacobi field $(J(\varphi))_{\varphi\in \D}$, which will have a simpler form.  

Let us consider the Hilbert space $\ell_2$ spanned by the
orthonormal basis
 $(e_n)_{n=0}^\infty$ with $$e_n=(0,\dots,0,\underbrace{1}_{(n+1)
 \text{-st place}},0,0\dots).$$ Consider the tensor product $\ell_2\otimes 
 L^2(X,\sigma)$, and let $$ {\cal F}(\ell_2\otimes 
 L^2(X,\sigma))=\bigoplus_{n=0}^\infty (\ell_2\otimes 
 L^2(X,\sigma))^{\hat\otimes n}n!$$ be the (usual) symmetric
Fock space over $\ell_2\otimes 
 L^2(X,\sigma)$.

Denote by $\ell_{2,0}$ the dense subset of $\ell_2$ consisting of
all finite vectors, i.e., $$
\ell_{2,0}{:=}\{(f^{(n)})_{n=0}^\infty: \exists N\in\Z_+\mathrm{\
such\ that\  } f^{(n)}=0\ \mathrm{for\ all\ }n\ge N\}.$$ The Jacobi
matrix $J$ determines a linear symmetric operator in $\ell_2$ with
domain $\ell_{2,0}$ by the following formula:
\begin{equation}\label{ifgvd} J e_n =  b_{n+1}e_{n+1}+ a_n
e_n+ b_n e_{n-1},\qquad n\in\Z_+,\ e_{-1}{:=}0.\end{equation} 
Denote by $J^+,J^0,J^-$ the corresponding creation, neutral, and annihilation operators in $\ell_{2,0}$, so that $J=J^++J^0+J^-$. 

Denote by $\Phi$ the linear subspace of ${\cal F}(\ell_2\otimes L^2(X,\sigma))$ that is the linear span of the vacuum vector
$(1,0,0,\dots)$
 and vectors of the form $(\xi\otimes \varphi)^{\otimes n}$, where $\xi\in \ell_{2,0}$, $\varphi\in\D$, and $n\in\N$. The set $\Phi$ is evidently a dense subset of ${\cal F}(\ell_2\otimes L^2(X,\sigma))$.

Now, for each $\varphi,\psi \in\D$ and $\xi\in\ell_{2,0}$, we set
\begin{align}
A^+(\varphi)(\xi\otimes\psi)^{\otimes n}&{:=}(e_0\otimes\varphi)\hat\otimes(\xi\otimes\psi)^{\otimes n}\notag\\&\quad+
n((J^+\xi)\otimes(\varphi\psi))\hat\otimes (\xi\otimes\varphi)^{\otimes(n-1)},\notag\\
A^0(\varphi)(\xi\otimes\psi)^{\otimes n}&{:=}
n((J^0\xi)\otimes(\varphi\psi))\hat\otimes (\xi\otimes\varphi)^{\otimes(n-1)},\notag\\
A^-(\varphi)(\xi\otimes\psi)^{\otimes n}&{:=}n\la\xi,e_0\ra\la\varphi,\psi\ra(\xi\otimes\varphi)^{\otimes(n-1)}
\notag \\&\quad+
n((J^-\xi)\otimes(\varphi\psi))\hat\otimes (\xi\otimes\varphi)^{\otimes(n-1)},\label{trd}
\end{align}
and then extend these operators by linearity  to the whole $\Phi$.
Thus, \begin{align*}
A^+(\varphi)&=a^+(e_0\otimes\varphi)+a^0(J^+\otimes\varphi),\\
A^0(\varphi)&=a^0(J^0\otimes\varphi),\\
A^-(\varphi)&=a^-(e_0\otimes\varphi)+a^0(J^-\otimes\varphi),
\end{align*}
where $a^+(\cdot)$, $a^0(\cdot)$, $a^-(\cdot)$ are the usual creation, neutral, and annihilation operators in $\mathcal F(\ell_2\otimes L^2(X,\sigma))$. (In fact, under, e.g., $a^0(J^+\otimes\varphi)$ we understand the differential second quantization of the operator $J^+\otimes \varphi$ in $\ell_2\otimes L^2(X,\sigma)$, which, in turn, is the tensor product of the operator $J^+$ in $\ell_2$ defined on $\ell_{2,0}$ and the operator of multiplication by $\varphi$ in $
L^2(X,\sigma)$ defined on $\D$.)
Note also that \begin{align*} A(\varphi){:=}&A^+(\varphi)+A^0(\varphi)+A^-(\varphi) \\
=&a^+(e_0\otimes\varphi)+a^0((J^++J^0+J^-)\otimes\varphi)+a^-(e_0\otimes\varphi)
\\=&a^+(e_0\otimes\varphi)+a^0(J\otimes\varphi)+a^-(e_0\otimes\varphi).\end{align*}

In the following theorem, for a linear operator $A$ in a Hilbert space $H$, we denote by $\overline A$ its closure (if it esxists). 

\begin{theorem}
There exists a unitary operator $$ I:{\frak F}\to {\cal F}(\ell_2\otimes 
 L^2(X,\sigma))$$ for which the following assertions hold\rom.
Let $J^+(\varphi)$\rom, $J^0(\varphi)$\rom, and $J^-(\varphi)$\rom, 
$\varphi\in\D$\rom, be linear operators in $\frak F$ with domain $\Ffin(\D)$ as in Theorem \rom{3.} Then\rom, for all $\varphi\in\D$\rom,
\begin{align*}
I\overline{J^+(\varphi)}I^{-1}&= \overline{A^+(\varphi)},\\
I\overline{J^0(\varphi)}I^{-1}&= \overline{A^0(\varphi)},\\
I\overline{J^-(\varphi)}I^{-1}&= \overline{A^-(\varphi)},
\end{align*} 
  and $$ I J(\varphi) I^{-1}= \overline{A(\varphi)}.$$

\end{theorem}

\begin{remark}{\rm
Note, however, that the image of ${\frak F}_n$ under $I$ does not coincide with the subspace  $(\ell_2\otimes L^2(X,\sigma))^{\hat\otimes n}n!$
of the Fock space ${\cal F}(\ell_2\otimes 
 L^2(X,\sigma))$.
}\end{remark}

The proof of Theorem 4 is a straightforward generalization of the proof of Theorem~1 in \cite{LSWN}, so we only outline it. 

First, we recall the classical unitary isomorphism between the usual 
Fock space over $L^2(\R\times X,\nu\otimes \sigma)$ and $L^2(\D',\mu)$:
$$ {\cal U}_1:{\cal F}(L^2(\R\times X,\nu\otimes\sigma))\to L^2(\D',\mu).$$ This isomorphism was broved by It\^o \cite{Ito} and extended in \cite{L} to a general L\'evy process on $X$. Note also that
$$ L^2(\R\times X,\nu\otimes\sigma)=L^2(\R,\nu)\otimes L^2(X,\sigma).$$ Next, we have the unitary operator $$ {\cal U}_2 :
L^2(\R,\tilde\nu)\to L^2(\R,\nu)$$ defined by $$ ({\cal U}_2 f)(s):=
\frac1s\,f(s). $$ 

Let $$ {\cal U}_3:\ell_2\to L^2(\R,\tilde\nu)$$ be the Fourier transform in generalized joint eigenvectors of the Jacobi matrix $J$, see \cite{b}. The ${\cal U}_3$ can be characterized 
as a unique unitary operator which maps the vector $(1,0,0,\dots)$
into the function identically equal to 1, and which maps the closure $\overline J$ of the symmetric operator $J$ in $\ell_2$ into the multiplication operator by the variable $s$. 

Let $$ {\cal U}_4:\ell_2\otimes L^2(X,\sigma)\to L^2(\R\times X,\nu\otimes \sigma)$$ be given by $${\cal U}_4:=({\cal U}_2{\cal U}_3)\otimes \operatorname{id},$$ where $\operatorname{id}$
denotes the identity operator. Using ${\cal U}_4$, we naturally construct the unitary operator $${\cal U}_5:{\cal F}(\ell_2\otimes L^2(X,\sigma))
\to {\cal F}(L^2(\R\times X,\nu\otimes \sigma)).$$

We  now define a unitary operator $$ I:= {\cal U}{\cal U}_1^{-1}
{\cal U}_5^{-1}: {\frak F}\to {\cal F}(\ell_2\otimes L^2(X,\sigma)).$$
Then, the assertions of Theorem 4 about the unitary operator $I$
follow from Theorem~3, the construction of the unitary operator $I$
(see in particular Theorem~3.1 in \cite{L}), \eqref{trd}, and a limiting procedure.

\end{document}